\documentclass[12pt]{article}
\begin{document}
\bibliographystyle{plain}

\input{psfig}
\newtheorem{theorem}{Theorem}[section]
\newtheorem{corollary}{Corollary}[section]
\newtheorem{note}{Note}[section]
\newtheorem{lemma}{Lemma}[section]
\newtheorem{definition}{Definition}[section]
\newtheorem{conjecture}{Conjecture}[section]
\newtheorem{proposition}{Proposition}[section]
\font\boldsets=msbm10
\def\mod{\mbox{ mod }}
\def\polylog{\mbox{polylog}}
\def\sign{\mbox{ sign}}
\def\Pr{\mbox{ Pr}}
\def\E{{\hbox{\boldsets  E}}}
\def\F{{\hbox{\boldsets  F}}}
\def\N{{\hbox{\boldsets  N}}}
\def\R{{\hbox{\boldsets  R}}}
\def\C{{\hbox{\boldsets  C}}}
\def\Z{{\hbox{\boldsets  Z}}}
\def\blackslug{\hbox{\hskip 1pt \vrule width 4pt height 8pt depth 1.5pt
\hskip 1pt}}

\newcommand {\half}{\frac 1 2}
\newcommand{\ignore}[1]{}

\newcommand{\bbox}{\vrule height7pt width4pt depth1pt}
\def\QED{\quad\blackslug\lower 8.5pt\null\par}
\def\proof{\par\penalty-1000\vskip .5 pt\noindent{\bf Proof\/: }}
\newcommand{\prb}{\begin{problem}}
\newcommand{\eprb}{\bbox\end{problem}}
\newcommand {\upint}[1]   {\lceil {#1} \rceil}

\title
{Constructing expander graphs by 2-lifts and discrepancy vs. spectral gap}
\author{
Yonatan Bilu and Nathan Linial
\thanks{
Institute of Computer Science,
Hebrew University Jerusalem 91904 Israel
\mbox{\em \{johnblue,nati\}@cs.huji.ac.il}. This research is supported by the Israeli Ministry of Science and the Israel Science Foundation.
}
}
\maketitle


\begin{abstract}
We present a new explicit construction for expander graphs with nearly
optimal spectral gap. The construction is based on a series of $2$-lift 
operations.
Let $G$ be a graph on $n$ vertices. A $2$-lift of $G$ is a graph
$H$ on $2n$ vertices, with a covering map
$\pi:H \rightarrow G$. It is not hard to see that all eigenvalues
of $G$ are also eigenvalues of $H$. In addition, $H$ has $n$ ``new'' 
eigenvalues. We conjecture that every $d$-regular graph has a $2$-lift
 such that all new eigenvalues are in the range 
$[-2\sqrt{d-1},2\sqrt{d-1}]$ (If true, this is tight
, e.g. by the Alon-Boppana bound). 
Here we show that every graph of maximal degree $d$ has a $2$-lift 
such that all ``new'' eigenvalues are in the range 
$[-c \sqrt{d \log^3d}, c \sqrt{d \log^3d}]$ for some constant $c$.
This leads to a polynomial time algorithm for constructing 
arbitrarily large
$d$-regular graphs, with second eigenvalue $O(\sqrt{d \log^3 d})$.\\
The proof uses the following lemma (Lemma \ref{wonder}): 
Let $A$ be a real symmetric
matrix with zeros on the diagonal. Let $d$ be such that the $l_1$
norm of each row in $A$ is at most $d$. 
Suppose that $\frac {|xAy|} {||x||||y||} \leq \alpha$
for every $x,y \in \{0,1\}^n$ with $<x,y>=0$.
Then the spectral radius of $A$ is $O(\alpha (\log(d/\alpha) + 1))$.
An interesting consequence of this lemma is a converse to the 
Expander Mixing Lemma.\\
\end{abstract}

{\em Keywords:} Lifts, Lifts of Graphs, Discrepancy, Expander Graphs,
Signed Graphs. 

\pagebreak

\section{Introduction}
An $d$-regular graph is called a {\em $\lambda$-expander},
if all its
eigenvalues but the first are in $[-\lambda,\lambda]$. Such graphs are 
interesting when $d$ is fixed, $\lambda < d$, and the number of vertices 
in the graph tends to infinity. Applications of such graphs in 
computer science and discrete mathematics are many, see for example 
\cite{course} for a survey.\\ 
It is known that random $d$-regular graphs are good expanders
(\cite{BS}, \cite{FKS}, \cite{Fri03}), yet many 
applications require an explicit construction. Some known construction 
appear in \cite{Mar73}, 
\cite{GaGa}, \cite{AlMi}, \cite{LPS}, \cite{AGM}, \cite{Mar88}, \cite{Ajtai} 
and \cite{RVW}).
The Alon-Boppana bound says that $\lambda \geq  2\sqrt{d-1}-o(1)$ 
(cf. \cite{Nili}). The graphs of \cite{LPS} and \cite{Mar88} satisfy
$\lambda \leq  2\sqrt{d-1}$, for infinitely many values of $d$, and are
constructed very efficiently. However, the analysis of the eigenvalues
in these construction relies on deep mathematical results.
Thus, it is 
interesting to look for construction whose analysis is elementary.\\
The first major step in this direction is a construction based on 
iterative use of the zig-zag product \cite{RVW}. 
This construction is simple to analyze,
and is very explicit, yet the eigenvalue bound falls somewhat short of what
might be hoped for. The graphs constructed
with the zig-zag product have second eigenvalue $O(d^{3/4})$, which
can be improved, with some
additional effort to $O(d^{2/3})$. Here we introduce an iterative 
construction based on 
$2$-lifts of graphs, which is close to being optimal and gives
$\lambda = O(\sqrt{d \log^3 d)}$.\\
A graph
$\hat{G}$ is called a $k$-lift of a ``base graph'' $G$ if there is a $k:1$
 {\em covering map} $\pi:V(\hat{G})\rightarrow V(G)$. Namely,
if $y_1,\ldots,y_d \in G$ 
are the neighbors of $x \in G$, 
then every $x' \in \pi^{-1}(x)$ has exactly one
vertex in each of the subsets $\pi^{-1}(y_i)$.
See \cite{AmLi02} for a general introduction to graph lifts. \\
The study of lifts of graphs has focused so far mainly on random lifts
\cite{AmLi02, AmLi, AmLiMa, LiRo, Fri}.
In particular, Amit and
Linial show in \cite{AmLi} that w.h.p. a random $k$-lift has a strictly
positive edge expansion. It is not hard to see that the eigenvalues of the 
base graph
are also eigenvalues of the lifted graph. These are called
by Joel Friedman the ``old'' eigenvalues of the lifted graph.
In \cite{Fri} he shows that
 w.h.p. a random $k$-lift of a $d$-regular graph on $n$ vertices
is ``weakly Ramanujan''. Namely, that all 
eigenvalues but, perhaps, those of the base graph, are, in absolute value,
$O(d^{3/4})$. In both cases the probability tends to $1$ as $k$ tends to
infinity.\\
Here we study $2$-lifts of graphs.
We conjecture that every $d$ regular graph has a $2$-lift with all new 
eigenvalues at most $2\sqrt{d-1}$ in absolute value.
It is not hard to show (e.g., using the Alon-Boppana bound~\cite{Nili})
that if this conjecture is true, it is tight.
We prove (in Theorem \ref{main-thm}) a slightly weaker result;
every graph of maximal degree
$d$ has a $2$-lift with all new eigenvalues 
$O(\sqrt{d \log^3d})$ in absolute value. 
Under some natural assumptions on the base graph, such a $2$-lift can be
found efficiently. This leads to a polynomial time algorithm for constructing 
families of $d$-regular expander graphs, with second eigenvalue 
$O(\sqrt{d \log^3d})$. 

A useful property of expander graphs is the so-called Expander Mixing Lemma.
Roughly, this lemma states that the number of edges between 
two subsets of vertices in an expander graph
is what is expected in a random graph, up to an
additive error that depends on the second eigenvalue.\\
A key lemma in this paper (Lemma \ref{wonder}) shows a close connection
between the combinatorial discrepancy in a symmetric martix, and its spectral
radius. This key lemma
implies the following converse to the Expander Mixing Lemma:
Let $G$ be a $d$-regular graph on $n$ vertices,
such that for every two subsets of vertices, $A$ and $B$, 
$|e(A,B) - d|A||B|/n| \leq \alpha \sqrt{|A||B|}$ for
some $\alpha < d$. Then all eigenvalues
of $G$ but the first are, in absolute value, $O(\alpha \log(d/\alpha))$. 
The fact that the bound is tight up to a logarithmic factor is 
surprising. It is known that expansion implies a spectral gap
(cf. \cite{Al86}), but the actual bounds are weak, and indeed
expansion alone does not imply strong bounds on the spectral gap
\cite{Kahale}.

The paper is organized as follows. After defining the basic objects -
expander graphs, signed graphs and $2$-lifts - in section \ref{sec-def}, 
we present the
main results in section \ref{sec-main}. In sub-section \ref{subsec-prel}
we observe that the spectrum of $2$-lifts has a simple characterization,
which suggests an iterative construction of expander graphs. 
It reduces the problem of constructing expander graphs to finding
a {\em signing} of the edges with a small spectral radius. 
In sub-section
\ref{subsec-exist} we show that such a signing always exists. In sub-section
\ref{subsec-derandom} we show how to find such a signing efficiently.
An alternative method is given in section \ref{sec-aki}, which leads to 
a somewhat stronger notion of explicitness.
Finally, in section \ref{sec-lme} we prove the converse to the Expander Mixing 
Lemma mentioned above.

\section{Definitions}\label{sec-def}
Let $G=(V,E)$ be a graph on $n$ vertices, and let $A$ be its adjacency
matrix. Let $\lambda_1 \geq \lambda_2 \geq \ldots \geq \lambda_n$ 
be the eigenvalues of $A$. We denote by 
$\lambda(G) = \max_{i=2,\ldots,n}|\lambda_i|$.
We say that $G$ is an $(n,d,\mu)-expander$ if $G$ is $d$-regular, and 
$\lambda(G) \leq \mu$. If $\lambda(G) \leq 2\sqrt{d-1}$ we say
that $G$ is {\em Ramanujan}. 

A {\em signing} of the edges of $G$ is a function  
$s:E(G) \rightarrow \{-1,1\}$. 
The {\em signed adjacency matrix} of a graph $G$ with a 
signing $s$ has rows and columns indexed
by the vertices of $G$. The $(x,y)$ entry
is $s(x,y)$ if  $(x,y) \in E$ and $0$ otherwise.\\
A {\em $2$-lift} of $G$, associated with a signing $s$, is a graph  $\hat{G}$
defined as follows. Associated with every
vertex $x \in V$ are two vertices, $x_0$ and $x_1$, called the {\em fiber}
of $x$. If
$(x,y) \in E$, and $s(x,y) = 1$ then the corresponding edges
in $\hat{G}$ are $(x_0,y_0)$ and  $(x_1,y_1)$. If
$s(x,y) = -1$, then the corresponding edges 
in $\hat{G}$ are $(x_0,y_1)$ and  $(x_1,y_0)$.
The graph $G$ is called the {\em base graph}, and  $\hat{G}$ a $2$-lift 
of $G$.
By the {\em spectral radius of a signing} we refer to the spectral radius
of the corresponding signed adjacency matrix. 
When the spectral radius of a signing of a $d$-regular graph is
$\tilde{O}(\sqrt{d})$ we say that the signing (or the lift) is 
{\em Quasi-Ramanujan}.

For $v,u \in \{-1,0,1\}^n$, denote $S(u) = supp(u)$, and 
$S(u,v) = supp(u) \cup supp(v)$. \\ 
It will be convenient to assume throughout that $V(G)=\{1,\ldots,n\}$.

\section{Quasi-Ramanujan $2$-Lifts and Quasi-Ramanujan Graphs}\label{sec-main}
\subsection{Preliminaries}\label{subsec-prel}
The eigenvalues of a $2$-lift 
of $G$ can be easily characterized in terms of the adjacency
matrix and the signed adjacency matrix:

\begin{lemma}\label{lem-1}
Let $A$ be the adjacency matrix of a graph $G$, and $A_s$ the signed adjacency 
matrix associated with a $2$-lift $\hat{G}$. Then every eigenvalue of
$A$ and every eigenvalue of $A_s$ are eigenvalues of $\hat{G}$. Furthermore,
the multiplicity of each eigenvalue of $\hat{G}$ is the sum of its 
multiplicities in $A$ and  $A_s$.
\end{lemma}
\proof  
It is not hard to see that 
the adjacency matrix of $\hat{G}$ is:
\begin{eqnarray*}
\hat{A} =
\left(
\begin{array}{cc}
A_1 & A_2 \\
A_2 & A_1
\end{array}
\right)
\end{eqnarray*}
Where $A_1$ is the adjacency matrix of $(V,s^{-1}(1))$ and 
$A_2$ the adjacency matrix of $(V,s^{-1}(-1))$. (So $A=A_1+A_2$, 
$A_s = A_1-A_2$).
Let $v$ be an eigenvector of $A$ with eigenvalue $\mu$. 
It is easy to check that $\hat{v} = (v\;v)$
is an eigenvector of $\hat{A}$ with eigenvalue $\mu$.\\
Similarly, if $u$ is an eigenvector of $A_s$ with eigenvalue $\lambda$,
then $\hat{u} = (u\;{-u})$ is an eigenvector of $\hat{A}$ with eigenvalue $\lambda.$\\
As the $\hat{v}$'s and $\hat{u}$'s are perpendicular and $2n$ in number, 
they are all the eigenvectors of $\hat{A}$.
\QED

We follow Friedman's (\cite{Fri}) nomenclature, and call the eigenvalues
of $A$ the {\em old} eigenvalues of $\hat{G}$,
and those of $A_s$ the {\em new} ones.

Consider the following scheme for constructing $(n,d,\lambda)$-expanders.
Start with $G_0=K_{d+1}$, the complete graph on $d+1$ vertices
\footnote{We could start with any small $d$-regular graph with a large
spectral gap. Such graphs are easy to find.}
. Its eigenvalues
are $d$, with multiplicity $1$, and $-1$, with multiplicity $d$.
We want to define
$G_i$ as a $2$-lift of $G_{i-1}$, such that all new eigenvalues are 
in the range $[-\lambda,\lambda]$. Assuming such a $2$-lifts always exist,
the $G_i$ constitute an infinite family of $(n,d,\lambda)$-expanders.\\
It is therefore natural to look for the smallest 
$\lambda = \lambda(d)$ such that
every graph of degree at most $d$ has a $2$-lift, with new eigenvalues
in the range $[-\lambda,\lambda]$. In other words, a signing with
spectral radius $\le \lambda$.

We note that $\lambda(d) \geq 2\sqrt{d-1}$ follows from
the Alon-Boppana bound. We next observe:
\begin{proposition}
Let $G$ be a $d$-regular graph which
contains a vertex that does not belong to any cycle of bounded length,
then no signing of $G$ has spectral radius below $2\sqrt{d-1}-o(1)$.
\end{proposition}
To see this, note first that all signing of a tree have the
same spectral radius. This follows e.g., from the easy fact
that any $2$-lift of a tree is a union of two disjoint trees,
isomorphic to the base graph. The assumption implies that $G$
contains an induced subgraph that is a full $d$-ary tree $T$ of
unbounded radius. The spectral radius of $T$ is $2\sqrt{d-1}-o(1)$.
The conclusion follows now from the interlacing principle of
eigenvalues.

\begin{figure}
\begin{center}
\begin{tabular}{c}
{\psfig{file=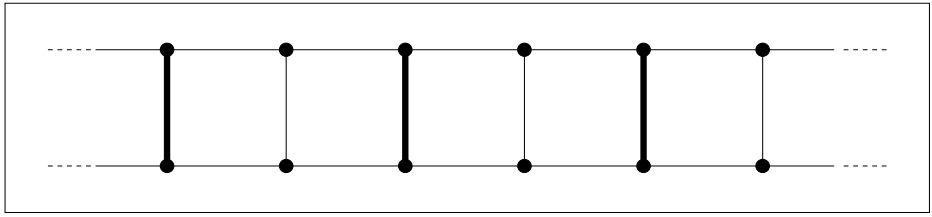,width=0.9\columnwidth}}
\end{tabular}
\end{center}
\caption{The Railway Graph. Edges where the signing is $-1$
are bold.}
\label{fig-1}
\end{figure}

There are several interesting examples of arbitrarily large $d$-regular 
graphs for which there is a signing with
spectral radius bounded away from $2\sqrt{d-1}$. One such example is
the $3$-regular graph $R$ defined as follows. 
$V(R) = \{0,\ldots,2k-1\} \times \{0,1\}$. For $i \in [2k], j \in \{0,1\}$, the
neighbors of $(i,j) \in R$ are $((i-1) \mod 2k,j)$, $((i+1) \mod 2k,j)$ and
$(i,1-j)$. 
Define $s$, a signing of $R$, to be $-1$ on the
edges $((2i,0),(2i,1))$, for $i \in \{0,\ldots,k-1\}$, and $1$ elsewhere
(see Figure \ref{fig-1}).
Let $A_s$ be the signed adjacency matrix. It is easy to see that 
$A_s^2$ is a matrix with $3$ on the diagonal, and two $1$'s in each
row and column. Thus, its spectral radius is $5$, and that
of $A_s$ is $\sqrt{5} < 2 \sqrt{2}$. 

\subsection{Quasi-ramanujan $2$-lifts for every graph}\label{subsec-exist}

We conjecture that every graph has a signing with small spectral radius:

\begin{conjecture} 
Every $d$-regular graph has a signing with spectral radius
at most $2\sqrt{d-1}$.
\end{conjecture}

We have numerically tested this conjecture quite extensively.\\
In this subsection we show a close upper bound:

\begin{theorem}\label{main-thm}
Every graph of maximal degree $d$
has a signing with spectral radius
$O(\sqrt{d \cdot \log^3d})$.
\end{theorem}

The theorem is an easy consequence of the following
two lemmata (along with Lemma \ref{lem-1}). The first shows
that with positive probability the Rayleigh quotient is small
for vectors in $v,u \in \{-1,0,1\}^n$. The second shows
how to conclude from this that all eigenvalues are small.

\begin{lemma}\label{lem-2}
For every graph of maximal degree $d$,
there exists a signing $s$ such that
for all $v,u \in \{-1,0,1\}^n$
the following holds:
\begin{eqnarray}\label{good}
\frac {|v^tA_su|} {||v|| ||u||} \leq 10 \sqrt{d \log d},
\end{eqnarray}
where $A_s$ is the signed adjacency matrix.
\end{lemma}

\begin{lemma}\label{wonder}
Let $A$ be an $n \times n$ real symmetric matrix such
that the $l_1$ norm of each row in $A$ is at most $d$.
Assume that for any two vectors, $u,v \in \{0,1\}^n$,
with $supp(u) \cap supp(v) = \emptyset$:
\[\frac {|uAv|} {||u|| ||v||} \leq \alpha,\]
and that all diagonal entries of $A$ are, in absolute value, 
$O(\alpha(log(d/\alpha)+1))$.
Then the spectral radius of $A$ is
$O(\alpha(log(d/\alpha)+1))$.
\end{lemma}

\proof (Lemma \ref{lem-2})
First note that it's enough to prove this for 
$u$'s and $v$'s such that the set $S(u,v)$ spans a connected subgraph.
Indeed, assume that the claim holds for all connected subgraphs and 
suppose that $S(u,v)$ is not connected. Split $u$ and
$v$ according to the connected components of $S(u,v)$, and apply the
claim to each component separately. 
Summing these up and using the Cauchy-Schwartz 
inequality,
we conclude that the claim for $u$ and $v$ as well. So henceforth we assume
that $S(u,v)$ is a connected.

Consider some $u, v \in  \{-1,0,1\}^n$. 
Suppose we choose the sign of each edge uniformly at 
random. Denote the resulting signed adjacency matrix by $A_s$,
and by $E_{u,v}$ the ``bad'' event that 
$\frac {|v^tA_su|} {||v|| ||u||} > 10 \sqrt{d \log d}$.
Assume w.l.o.g. that $|S(u)| \geq \half |S(u,v)|$. 
By the Chernoff inequality ($v^tA_su$ is the sum of independent variables,
attaining values of either ${\pm 1}$ or  ${\pm 2}$):
\begin{eqnarray*}\label{chernoff}
Pr[E_{u,v}] \leq 
&&2 exp(-\frac {100 d \log d |S(u)| |S(v)|} {8 |e(S(u),S(v))|}) \\
&&\leq 2 exp(-\frac {100 d \log d |S(u)| |S(v)|} {8 d |S(v)|}) \\
&&< d^{(-10|S(u,v)|)} 
\end{eqnarray*} 

We want to use the Lov\'asz
Local Lemma \cite{lll}, with the following dependency graph on the
$E_{u,v}$: There is an edge between $E_{u,v}$ and $E_{u',v'}$ iff
$S(u,v) \cap S(u',v') \neq \emptyset$. Denote $k = |S(u,v)|$. How many neighbors,
$E_{u',v'}$, does $E_{u,v}$ have, with $|S(u',v')| = l$? 

Since we are interested
only in connected subsets, this is clearly bounded by the number of rooted directed 
subtrees on $l$ vertices, with a root in $S(u,v)$.
It is known (cf. \cite{Knuth}) that 
there are at most $k {d(l-1) \choose l-1} \approx k d^{l-1}$ such trees
(a similar argument appears in \cite{FrMo}. The bound on the number of
trees is essentially tight by \cite{Al91}).

In order to apply the Local Lemma, 
we need to define for such $u$ and $v$ a numbers $0 \leq X_{u,v} < 1$.
It is required that:
\begin{eqnarray}\label{cond1}
X_{u,v} \prod_{(u',v'): E_{u,v} \sim E_{u',v'}}(1-X_{u',v'}) \geq 
d^{-10|S(u,v)|}.
\end{eqnarray}
Observe that for $S \subseteq [n]$ there are at most
$2^{4|S|}$ distinct pairs $v,u \in \{-1,0,1\}^n$ with 
such that $S(u,v) = S$.\\ 
For all $u,v$ set $X_{u,v} = d^{-3k}$,
where $k=|S(u,v)|$.
Then in (\ref{cond1}) we get:
\begin{eqnarray*}
&& X_{u,v} \cdot \prod_{(u',v'): E_{u,v} \sim E_{u',v'}}(1-X_{u',v'}) =
d^{-3k}\prod_{l=1}^n(1-d^{-3l})^{kd^l2^{4l}} \approx \\
&&d^{-3k}exp(-k\sum_{l=1}^nd^{-3l}d^{l}2^{4l}) \geq
d^{-3k}e^{-2k} > d^{-10k}
\end{eqnarray*}
as required. 
\QED

\proof (Lemma \ref{wonder})
For simplicity, assume that all diagonal entries of $A$ are zeros. We explain
at the end of the proof how to deal with a general matrix (but in fact
it does not matter for the purpose of this paper).

First note that our assumptions imply that for any  $u \in \{0,1\}^n$,
\[\frac {|uAu|} {||u||^2} \leq 2 \alpha:\] 
For any $u_1, u_2 \in \{0,1\}^n$ such $S(u_1) \cap S(u_2) = \emptyset$
we have that 
\begin{eqnarray}\label{asump}
|u_1Au_2| \leq \alpha ||u_1|| ||u_2||.
\end{eqnarray}
Let $u \in \{0,1\}^n$, and denote $k=|S(u)|$. Set $K = {k \choose k/2}$.
Summing up inequality (\ref{asump}) over all subsets of $S(u)$ of size
$k/2$, we have that:
\[\sum_{u_1:S(u_1) \subset S(u),|S(u_1)|=k/2}u_1 A (u-u_1) \leq K \alpha k/2.\]
For each $i \neq j \in S(u)$, $a_{i,j}$ is added up ${k-2 \choose k/2-1}$ 
times in the sum on the LHS, hence (since diagonal entries are zero):
\[ {k-2 \choose k/2-1} uAu \leq K \alpha k/2,\]
or: 
\[ uAu \leq \alpha 2k.\] 

Next, it follows
that for any  $u,v \in \{-1,0,1\}^n$, 
such that $S(u) = S(v)$, or $S(u) \cap S(v) = \emptyset$:
\[\frac {|uAv|} {||u|| ||v||} \leq 4 \alpha.\]
Fix $u,v \in \{-1,0,1\}^n$.
Denote $u=u^+-u^-$ and $v=v^+-v^-$, where 
$u^+, u^-, v^+, v^- \in \{0,1\}^n$, and $S(u^+)\cap S(u^-) = 
S(v^+) \cap (v^-) = \emptyset$.
\begin{eqnarray*}
&&|uAv| = |(u^+-u^-)A(v^+-v^-)|\\
&&\leq 2 \alpha (||u^+||||v^+|| + ||u^+||||v^-|| + ||u^-||||v^+|| +
||u^-||||v^-||)\\
&&\leq 4 \alpha
\sqrt{||u^+||^2||v^+||^2+||u^+||^2||v^-||^2+||u^-||^2||v^+||^2+||u^-||^2||v^-||^2}\\ 
&& = 4 \alpha
\sqrt{(||u^+||^2 + ||u^-||^2)(||v^+||^2+||v^-||^2)} = 4 \alpha ||u|| ||v||\\
\end{eqnarray*}
The first inequality follows from our assumption on vectors in  $\{0,1\}^n$, 
and the second from the $l_2$ to $l_1$ norm ratio.

Fix  $x \in \R^n$.
We need to show that 
$\frac {|xAx|} {||x||^2} = O(\alpha \log(d/\alpha))$. 
By losing only a multiplicative factor of $2$, we may assume that the
absolute value of every
non-zero entry in $x$ is a negative powers of $2$: 
Clearly we may assume that $||x||_{\infty} < \half$.
To bound the effect of rounding the coordinates, denote 
$x_i = \pm (1+\delta_i)2^{t_i}$, with $0 \leq \delta_i \leq 1$ and $t_i < -1$,
an integer. Now round $x$ to a vector $x'$ by choosing the value of $x'_i$ to be
$\sign(x_i) \cdot 2^{t_i+1}$ with probability $\delta_i$ and 
$\sign(x_i) \cdot 2^{t_i}$ with probability 
$1 - \delta_i$. The expectation of $x'_i$ is $x_i$. As the coordinates
of $x'$ are chosen independently, and the diagonal entries of $A$ are $0$'s,
the expectation of $x'Ax'$ is $xAx$. Thus, there's a rounding, $x'$,
of $x$, such that $|xAx| \leq |x'Ax'|$. Clearly $||x'||^2 \leq 2 ||x||^2$, so
$\frac {|xAx|} {||x||^2} \leq 2 \frac {|x'Ax'|} {||x'||^2}$.
 
Denote $S_i=\{j:x_j = \pm 2^{-i}\}$, $s_i=|S_i|$.
Denote by $k$ the maximal index $i$ such that $s_i > 0$.
Denote by $x^i$ the sign vector of $x$ restricted to $S_i$, that is,
the vector whose $j$'th coordinate is the sign of $x_j$ if $j \in S_i$,
and zero otherwise. By our assumptions, 
for all $1 \leq i \leq j \leq k$:
\begin{eqnarray}\label{alpha-cond}
|x^iAx^j| \leq \alpha \sqrt{s_i s_j}.
\end{eqnarray}
Also, since the $l_1$ norm of each row is at most $d$, for all 
$1 \leq i \leq k$:
\begin{eqnarray}\label{d-cond}
\sum_j|x^iAx^j| \leq d s_i.
\end{eqnarray}
We wish to bound: 
\begin{eqnarray}\label{ray-q}
\frac {|xAx|}  {||x||^2} \leq \frac {\sum_{i,j=1}^k |x^iAx^j| 2^{-(i+j)}} 
{\sum_i 2^{-2i} s_i}.
\end{eqnarray}
Denote $\gamma = \log(d/\alpha)$, $q_i = s_i 2^{-2i}$ and
$ Q = \sum_i q_i$. 
Add up inequalities (\ref{alpha-cond}) and (\ref{d-cond}) as follows. For $i=j$
multiply inequality (\ref{alpha-cond}) by $2^{-2i}$. When
 $i < j \leq i + \gamma$ multiply it by $2^{-(i+j)+1}$. Multiply
inequality (\ref{d-cond}) by $2^{-(2i+\gamma)}$. (We ignore inequalities
(\ref{alpha-cond}) when $j > i + \gamma$.)

We get that:
\begin{eqnarray*} 
&&\sum_i2^{-2i}|x^iAx^i| +
\sum_i\sum_{i < j \leq i+ \gamma}2^{-(i+j)+1}|x^iAx^j| + 
\sum_i 2^{-(2i+\gamma)} \sum_j|x^iAx^j|\\
&& \leq \sum_i \alpha q_i + 
\sum_i \sum_{i < j < i+ \gamma}2 \alpha \sqrt{q_i q_j} + 
\sum_i 2^{-\gamma}d \cdot q_i \\
&& \leq \alpha \sum_i  q_i +
\alpha\sum_i \sum_{i < j < i+ \gamma}(q_i + q_j) +
2^{-\gamma} d \sum_i  \cdot q_i \\
&& < (2^{-\gamma}d + 2\gamma\alpha) \sum_i q_i = 
(\alpha + \alpha \log(d/\alpha))Q.
\end{eqnarray*}
Note that the denominator in (\ref{ray-q}) is $Q$, so to prove the
lemma it's enough to show that the numerator,
\begin{eqnarray}\label{numer}
\sum_{i<j} 2^{-(i+j+1)}|x^iAx^j|  + \sum_i  2^{-2i}|x^iAx^i|,
\end{eqnarray}
is bounded by
\begin{eqnarray}\label{ineq-sum}
\sum_i2^{-2i}|x^iAx^i| +
\sum_i\sum_{i\leq j \leq i+ \gamma}2^{-(i+j)+1}|x^iAx^j| + 
\sum_i 2^{-(2i+\gamma)} \sum_j|x^iAx^j|.
\end{eqnarray}
Indeed, let us compare the coefficients of the terms $|x^iAx^j|$
in both expressions (Since $|x^iAx^j| = |x^jAx^i|$, it's enough
to consider $i \leq j$). 
For $i=j$ this coefficient 
is $2^{-2i}$ in (\ref{numer}), 
and $2^{-2i}+2^{-(2i+\gamma)}$ in (\ref{ineq-sum}).
For $i < j \leq i + \gamma$, it is $2^{-(i+j)+1}$
in (\ref{numer}), and
in (\ref{ineq-sum}) it is  $2^{-(i+j)+1}+2^{-(2i+\gamma)}+2^{-(2j+\gamma)}$.
For $j > i + \gamma$, in (\ref{numer}) the coefficient is again $2^{-(i+j)+1}$.
In (\ref{ineq-sum}) it is:
\[2^{-(2i+\gamma)}+2^{-(2j+\gamma)} > 2^{-(2i+\gamma)} \geq 2^{-(i+j) + 1}.\]

It remains to show that the lemma holds when the diagonal entries of $A$
are not zeros, but 
$O(\alpha(log(d/\alpha)+1))$ in absolute value.
Denote $B = A - D$,
with $D$ the matrix having the entries of $A$ on the diagonal, and zero
elsewhere. 
We have that for any two vectors, $u,v \in \{0,1\}^n$,
with $supp(u) \cap supp(v) = \emptyset$:
\[\frac {|uAv|} {||u|| ||v||} \leq \alpha.\]
For such vectors, $uBv = uAv$, so by applying the lemma to $B$,
we get that its spectral radius is
$O(\alpha(log(d/\alpha)+1))$.
By the assumption on the diagonal entries of $A$, the spectral radius of $D$
is
$O(\alpha(log(d/\alpha)+1))$ as well. The spectral radius of $B$ is at most
the sum of these bounds - also 
$O(\alpha(log(d/\alpha)+1))$.
\QED

Note that Lemma \ref{wonder} is tight up to constant factors. 
To see this, consider the 
$n$-dimensional
vector $x$ whose $i$'th entry is $1/\sqrt{i}$. Let $A$ be the outer
product of $x$ with itself, that is, the matrix whose $(i,j)$'th 
entry is $1/\sqrt{i \cdot j}$. Clearly $x$ is an eigenvector of $A$
corresponding to the eigenvalue $||x||^2 = \Theta(\log(n))$. 
Also, the sum of each row in $A$ is $O(\sqrt{n})$. To prove
that the lemma is essentially tight, we need to show that 
$max_{u,v \in \{0,1\}^n} \frac {uAv} {||u||||v||}$ is constant.
Indeed, fix $k,l \in [n]$.
Let $u,v \in \{0,1\}^n$ be such that $||u|| = k$ and $||v||=l$. 
As the entries of
$A$ are decreasing along the rows and the columns, 
$uAv$ is maximized for such vectors when their support is the first
$k$ and $l$ coordinates. For these optimal vectors, 
$uAv = \Theta(\sqrt{k\cdot l})$. Thus, 
\[max_{u,v \in \{0,1\}^n} \frac {uAv} {||u||||v||} = \Theta(1).\]
A more elaborate example appears in section \ref{sec-lme}.

\subsection{An explicit construction of quasi-ramanujan graphs}\label{subsec-derandom}

For the purpose of constructing expanders, it is enough to prove a weaker
version of Theorem \ref{main-thm}. Roughly, that every {\em expander} graph
has a $2$-lift with small spectral radius. In this sub-section
we show that when the base graph is a good expander (in the sense of
the definition below), then w.h.p. a 
random $2$-lift has a small spectral radius. We then derandomize
the construction to get a deterministic
polynomial time
algorithm for constructing arbitrarily large expander graphs.

\begin{definition}
We say that a graph $G$ on $n$ vertices is $(\beta,t)$-sparse 
if for every $u,v \in \{0,1\}^n$, with $|S(u,v)| \leq t$, 
\[uAv \leq \beta ||u|| ||v||.\]
\end{definition}

\begin{lemma}\label{lem-2e}
Let $A$ be the adjacency matrix of a $d$-regular 
$(\gamma(d),\log n)$-sparse $G$ graph on $n$ vertices, where
$\gamma(d) = 10 \sqrt{d \log d}$.
Then for a random signing of $G$ (where the sign of each edge is chosen
uniformly at random) the following hold w.h.p.:
\begin{enumerate}
\item\label{req-1}$\forall u,v \in \{-1,0,1\}^n: 
|uA_sv| \leq \gamma(d) ||u|| ||v||$.
\item\label{req-2}$\hat{G}$ is $(\gamma(d),1 +\log n)$-sparse
\end{enumerate}
where $A_s$ is the random signed adjacency matrix, and $\hat{G}$ is the
corresponding $2$-lift.
\end{lemma} 
\proof
Following the same arguments 
and notations as in the proof of Lemma \ref{lem-2}, we have that
there are at most $n \cdot d^k$ connected subsets of size $k$. With
probability at most $d^{-10k}$ requirement (\ref{req-1}) is violated
for a given pair $u,v$ such that $|S(u,v)| = k$. Since for each $S$ 
there are at most $2^{4|S|}$ pairs $u,v$ such that $S(u,v) = S$, by the union
bound, w.h.p. no pair $u,v$ such that $|S(u,v)| > \log n$ violates
(\ref{req-1}). If
$|S(u,v)| \leq \log n$ then by (\ref{req-2}) there are simply not 
enough edges between $S(u)$ and $S(v)$ for (\ref{req-1}) to be violated. 

Next, we show that w.h.p. (\ref{req-2}) holds as well.
Let $s$ be a signing, and define $A_1$, $A_2$ and $\hat{A}$ as in Lemma \ref{lem-1}.
Given $u=(u_1\; u_2),v =(v_1\; v_2)\in \{0,1\}^n \times \{0,1\}^n$,
we wish to prove that $uAv \leq \gamma(d) ||u|| ||v||$.
As in the proof of Lemma \ref{lem-2}
we may assume that $S(u,v)$ is connected - in fact, that it is
connected via the edges between $S(u)$ and $S(v)$. Hence, we may assume
that the ratio of the sizes of these subsets is at most $d$.
Define
$x = u_1 \vee u_2$, $y = v_1 \vee v_2$, $x' = u_1 \wedge u_2$, and $y' = v_1 \wedge v_2$ (the characteristic vectors of $S(u_1,u_2)$, $S(v_1,v_2)$, 
$S(u_1)\cap S(u_2)$ and $S(v_1)\cap S(v_2)$).
It is not hard to verify that:
\begin{eqnarray}\label{ineq1}
u\hat{A}v = u_1A_1v_1 + u_1A_2v_2 + u_2A_2v_1 + u_2A_1v_2 \leq
xAy + x'Ay'.
\end{eqnarray}
If $|S(x,y)| \leq \log n$, then clearly $|S(x',y')| \leq \log n$
and from the assumption that
$G$ is $(\gamma(d),\log n)$-sparse
\[xAy + x'Ay' \leq \gamma(d)(\sqrt{|S(x)| |S(y)|} + \sqrt{|S(x')| |S(y')|}).\] 
Observe that $|S(u)| = |S(x)| + |S(x')|$ and $|S(v)| = |S(y)| + |S(y')|$,
so in particular $u\hat{A}v \leq \gamma(d)\sqrt{|S(u)| |S(v)|}$, and
requirement (\ref{req-2}) holds.\\
So assume 
$|S(x,y)| =  |S(u,v)| = \log n +1$. 
It is not hard to see that this entails 
$S(u_1,v_1) \cap S(u_2,v_2) = \emptyset$.
In other words, $S(u,v)$ contains at most one vertex from each fiber.
Hence, $x' = y' = \vec{0}$ and
$|S(u)| = |S(x)|$, $|S(v)| = |S(y)|$.\\
Denote $S = S(x,y)$, and assume w.l.o.g. that $|S(y)| > \half \log n$.
From  (\ref{ineq1}) $u\hat{A}v \leq xAy$, so it's enough to show that
$xAy \leq \gamma(d)\sqrt{|S(x)||S(y)|}$. If this is not the case, we can bound
the ratio between $|S(x)|$ and $|S(y)|$:
Since the graph is of maximal degree $d$ we have $xAy \leq d|S(X)|$. Hence,
$\frac{|S(x)|} {|S(y)|} > \frac {\gamma(d)^2} {d^2} = 
\frac {100 \log d} {d}$.\\ 
Observe that the 
edges between $S(u)$ and $S(v)$ in $\hat{G}$
originate from edges between $S(x)$ and $S(y)$ in $G$ in the following
way - 
for each edge between $S(x)$ and $S(y)$ in $G$ there is, with probability
$\half$, an edge between $S(u)$ and $S(v)$ in  $\hat{G}$.\\
Next we bound $xAy$.
Averaging over all $S \backslash \{i\}$, for $i \in S(y)$ we have that:
\[(|S(y)|-2)xAy \leq |S(y)|\gamma(d)\sqrt{(|S(y)|-1)|S(x)|}.\]
Hence the expectation of $u\hat{A}v$ is at most $\half c \gamma(d)||u|| ||v||$,
where 
\[c = \frac {\sqrt{|S(y)|(|S(y)|-1)}} {|S(y)|-2} \leq 1.1\] (assuming
$n$ is not very small). By the Chernoff bound, the probability that
$u\hat{A}v > \gamma(d) ||u|| ||v||$ is at most:
\begin{eqnarray*}
&&2exp(-\frac {0.9} {2.2} \gamma(d) ||u|| ||v||) \leq \\
&& exp(-0.2 \gamma(d) \frac {(\log n + 1)(10 \sqrt{\log d})} {\sqrt{d}}) = \\
&& exp(-20 \log d (\log n + 1)),
\end{eqnarray*}
Since $\frac{|S(x)|} {|S(y)|} > \frac {100 \log d} {d}$, 
and $|S(v)| > \half \log n$.\\
There are at most $d^{\log n+1} 4^{\log n+1}$ pairs $u,v$ with $S(u,v)$ connected and of size
$\log n+1$, so by the union bound, w.h.p., requirement (\ref{req-2}) holds.
\QED

\begin{corollary}\label{cor-2e}
Let $A$ be the adjacency matrix of a $d$-regular 
$(\gamma(d),\log n)$-sparse $G$ graph on $n$ vertices, where
$\gamma(d) = 10 \sqrt{d \log d}$.
Then there is a deterministic polynomial time algorithm for finding a
signing $s$ of $G$ such that the following hold:
\begin{enumerate}
\item The spectral radius of $A_s$ is $O(\sqrt{d \log^3 d})$.
\item $\hat{G}$ is $(\gamma(d),1 +\log n)$-sparse,
\end{enumerate}
where $A_s$ is the signed adjacency matrix, and $\hat{G}$ is the
corresponding $2$-lift.
\end{corollary} 
\proof 
Consider a random signing $s$.
For each closed path $p$ in $G$ of length $l = 2 \upint{\log n}$ 
define a random variable
$Y_p$ equal to the product of the signs of its edges. 
From lemmas \ref{lem-2e} and \ref{wonder}, the expected value of the 
trace of $A_s^l$, which is the expected value of the sum of these
variables, is $O(\sqrt{d \log^3 d})^l$ (since $l$ is even the sum is always
positive). For each $u,v \in \{0,1\}^n$, with
$|S(u,v)|=\log n + 1$, and $S(u,v)$ connected, define $Z_{u,v}$ to be
$d^l$ if $u\hat{A}v \leq \gamma(d) ||u|| ||v||$, and $0$ otherwise. 
In the proof of Lemma \ref{lem-2e} we've seen that the probability that
$Z_{u,v}$ is not $0$ is at most $d^{-10\log n}$, thus the expected value
of $Z_{u,v}$ is at most $d^{-8\log n}$. Let $Z$ be the sum of the $Z_{u,v}$'s.
Recall that there are at most $n (4d)^{\log n + 1}$ pairs $(u,v)$ such that 
$|S(u,v)|=\log n + 1$, and $S(u,v)$ connected. Hence, the expected value
of $Z$ is less than $d^{-6\log n}$.

Let $X = Y + Z$.
Note that the expected value of $X$ is approximately
that of $Y$, $O(\sqrt{d \log^3 d})^l$.
The expectation of $Y_p$ and $Z_{u,v}$ can be easily computed even when
the sign of some of the edges is fixed, and that of the other is chosen at
random. As there is only a polynomial number of variables, using the method
of conditional probabilities (cf. \cite{AlSp})
one can find a signing $s$ such that the value
of $X$ is at most its expectation.
For this value of $X$,
$tr(A_s^l) = Y = O(\sqrt{d \log^3 d})^l$, and
$Z = 0$ since if $Z \neq 0$ then $Z \geq d^l$. 
Clearly, the spectral radius of $A_s$ is $O(\sqrt{d \log^3 d})$. 
In the proof of Lemma 
\ref{lem-2e} we've seen that if $G$ is $(\gamma(d),\log n)$-sparse then so 
is $\hat{G}$, for any signing of $G$. For our choice of $s$ all $Z_{u,v} = 0$,
hence $\hat{G}$ is actually $(\gamma(d),\log n + 1)$-sparse. 
\QED

An alternative method for derandomization, using an almost $k$-wise 
independent sample space, is given in the appendix.

Recall the construction from the beginning of this section. 
Start with a $d$-regular
graph $G_0$ which is an $(n_0,d,\mu)-expander$ expander, 
for $\mu  = 10 \sqrt{d \log d}$ and $n_0 > d \log^2 n_0$.
From the Expander Mixing Lemma (cf. \cite{AlSp}), 
$G_0$ is $(\mu, \log n_0)$-sparse. Iteratively chose $G_{i+1}$
to be a $2$-lift
of $G_i$ according to Corollary \ref{cor-2e}, for $i=1,\ldots,\log(n/n_0)$.
Clearly this is a polynomial time algorithm that yields a
$(n,d,O(\sqrt{d \log^3 d}))$-expander graph.

\subsection{Random $2$-lifts}
Theorem \ref{main-thm} states that for every graph there exists a signing
such that the spectral radius of the signed matrix is small. The proof shows
that for a random signing, this happens with positive, yet exponentially small,
probability. The following example shows the limitations of this argument
and in particular, that there exist graphs for which a random signing
almost surely fails to give a small spectral radius.

Consider a graph composed of $n/(d+1)$ disjoint copies of $K_{d+1}$ (the
complete graph on $d+1$ vertices).
If all edges in one of the components are equally signed,
then $A_s$ has spectral radius $d$. For $d$ fixed and $n$ large, this
event will occur with high probability. Note that connectivity is not
the issue here - it is easy to modify this
example and get a connected graph for which, w.h.p., the spectral radius
of $A_s$ is $\Omega(d)$.

However, for a random $d$-regular graph, it is true that a random $2$-lift
will, w.h.p., yield a signed matrix with small spectral radius. This follows
from the fact that, w.h.p., a random $d$-regular graph is an 
$(n,d,O(\sqrt{d}))$-expander (\cite{FKS,Fri91,Fri03}). 
In particular, by the Expander Mixing Lemma,
it is $(O(\sqrt{d}),\log n)$-sparse. By Lemma \ref{lem-2e}, w.h.p., a random
$2$-lift yields a signed matrix with small spectral radius.

\section{A stronger notion of explicitness}\label{sec-aki}
In this section we suggest an alternative derandomization scheme
of section \ref{subsec-derandom}.
We use the construction of Naor and Naor \cite{Naor} of a small, almost 
$k$-wise independent sample space. This derandomization scheme leads to
a construction which, in a sense, is more explicit than that in section
 \ref{subsec-derandom}.
\subsection{Derandomization using an almost $k$-wise independent sample space}
Recall that in the proof of Corollary \ref{cor-2e} we defined 
two types of random variables: Let $l = 2 \upint{\log n}$. For each closed
path $p$ of length $l$, $Y_p$ is the product of the signs of the edges of $p$.
For each $u,v \in \{0,1\}^n$, with
$|S(u,v)|=\log n + 1$, and $S(u,v)$ connected, let $Z_{u,v}$ be
$d^l$ if $u\hat{A}v \leq \gamma(d) ||u|| ||v||$, and $0$ otherwise.
Define $X$ to be the sum of all these random variables.\\
For brevity it will be convenient to make the following ad-hoc definitions:
\begin{definition}
A signing $s$ of a $d$-regular graph $G$ is {\em $(n,d)$-good}, if
the spectral radius of $A_s$ is $O(\sqrt{d \log^3 d})$ and 
$\hat{G}$ is $(\gamma(d),1 +\log n)$-sparse.\\
A $d$-regular graph $G$ is a {\em $(n,d)$-good expander}, 
if it's an
$(n,d,O(\sqrt{d \log^3 d})$-expander, and is $(\gamma(d),1 +\log n)$-sparse.
\end{definition}
The proof showed
that finding a good signing is equivalent to finding a 
signing such that the value of $X$ is at most its expected value.
We now show that this conclusion is also true when rather than choosing the
sign of each edge uniformly and independently, we choose the signing from
an $(\epsilon,k)$-wise independent sample space, with $k = d \log n$ and
$\epsilon = d^{- 2d \log n}$.
\begin{definition} (\cite{Naor})
Let $\Omega_m$ be a sample space of $m$-bit strings, and let 
$S = s_1 \ldots s_m$ be chosen uniformly at random from $\Omega_m$.
We say that $\Omega_m$ is an $(\epsilon,k)$-wise independent sample space
if for any $k' \leq k$ positions $i_1 < i_2 \ldots < i_{k'}$,
\[\sum_{\alpha \in \{-1,1\}^{k'}}|\Pr[ s_{i_1} \ldots s_{i_{k'}} = \alpha] 
- 2^{-k'}| < \epsilon.\] 
\end{definition}
Naor and Naor \cite{Naor} suggest an explicit construction of such sample
spaces. When $k=O(\log m)$ and $1/\epsilon = poly(m)$, the size of 
the sample space is polynomial in $m$ (simpler constructions are also given in
\cite{AGHP}).\\
We shall immediately see that the expected value of $X$ does not change 
significantly when the signing is chosen from such a sample space. Hence, an 
alternative way of efficiently finding a good signing is to go over the entire
sample space. For at least one point in it, the value of $X$ is at most its 
expected value, and thus the signing is good. 
\begin{lemma}
Let $m=dn/2$, $k = d \log n$ and $\epsilon = d^{- 2d \log n}$. Let
$\Omega_m$ be an $(\epsilon,k)$-wise independent sample space. Let $X$
be as in the proof of Corollary \ref{cor-2e}. Let $U_m$ be the 
uniform distribution on $m$ bits. Then
\[|\E_{\Omega_m}[X] - \E_{U_m}[X]| = o(1).\] 
\end{lemma}
\proof
Recall that $X = \sum_p Y_p + \sum_{u,v} Z_{u,v}$, 
where $Y_p$ and $Z_{u,v}$ are as
above, the first sum is over all closed paths of length $l = 2 \upint{\log n}$,
and the second sum is over all $u,v \in \{0,1\}^n$, with $|S(u,v)|=\log n + 1$
and $S(u,v)$ connected. Hence
\[|\E_{\Omega_m}[X] - \E_{U_m}[X]| \leq \sum_p |\E_{\Omega_m}[Y_p] - 
\E_{U_m}[Y_p]| + \sum_{u,v} |\E_{\Omega_m}[Z_{u,v}] - 
\E_{U_m}[Z_{u,v}]|.\]
Let $p$ be a path of length $l$, and denote the edges that appear in it an odd 
number of times by $i_1, \ldots i_{l'}$, for some $l' < l$. 
Let $s_{i_1}, \ldots s_{i_{l'}}$ be the signs of these edges. Then the
value of $Y_p$ is $\prod_{j=1}^{l'} s_{i_j}$, and (for every distribution)
\[\E[Y_p] = \sum_{\alpha \in  \{-1,1\}^{l'}} 
\Pr[s_{i_1}, \ldots s_{i_{l'}} = \alpha] \cdot \prod_{j=1}^{l'} \alpha_j.\]
Thus,  
\begin{eqnarray*}
&&|\E_{\Omega_m}[Y_p] - \E_{U_m}[Y_p]| =\\ 
&&|\sum_{\alpha \in  \{-1,1\}^{l'}} (\prod_{j=1}^{l'} \alpha_j)
(\Pr_{\Omega_m}[s_{i_1}, \ldots s_{i_{l'}} = \alpha] -
 \Pr_{U_m}[s_{i_1}, \ldots s_{i_{l'}} = \alpha])| \leq\\
&&|\sum_{\alpha \in  \{-1,1\}^{l'}} 
(\Pr_{\Omega_m}[s_{i_1}, \ldots s_{i_{l'}} = \alpha] - 2^{-l'})| < \epsilon.
\end{eqnarray*}  
As there are less than $d^l$ closed paths $p$ of length $l$, 
$\sum_p |\E_{\Omega_m}[Y_p] - \E_{U_m}[Y_p]| < \epsilon d^l = o(1)$.
A similar argument shows that 
$|\sum_{u,v} |\E_{\Omega_m}[Z_{u,v}] - \E_{U_m}[Z_{u,v}]| = o(1)$ as well.
\QED

In fact, it follows that w.h.p. (say, $1-\frac 1 {n^2}$ for an appropriate 
choice of $\epsilon$), choosing an element uniformly at random from 
$\Omega_m$ leads to a good signing.

\subsection{A probabilistic strongly explicit construction}
The constructions of section \ref{subsec-derandom} and of the previous 
subsection are explicit in the sense that given $n$ and $d$ they suggest
a polynomial (in $n$) time algorithm for constructing an $(n,d)$-good 
expander.
However, in some applications a stronger notion of explicitness is required.
Namely, an algorithm that given $n$, $d$ and $i,j \in [n]$, decides
in time $\polylog(n)$ whether $i$ and $j$ are adjacent. We do not know how
to achieve such explicitness using the $2$-lifts schema, but we can do so
probabilistically. Consider the following:
\begin{definition}
Let $f_n:\{0,1\}^t \times [n] \times [n] \rightarrow \{0,1\}$, with
$t = O(\log n)$. Given $r \in \{0,1\}^t$, $f_n$ defines a graph $G_{f_n}(r)$,
on $n$ vertices,
where $i$ and $j$ are adjacent iff $f_n(r,i,j) = 1$. We say that $f_n$ is
a {\em $\delta$-probabilistic strongly explicit} description of an 
$(n,d)$-good expander graph,
if given $n$, $f_n$ can be computed in time $\polylog(n)$, and, with
probability at least $1-\delta$ (over
a uniform choice of $r$), $G_{f_n}(r)$ is an $(n,d)$-good expander graph.  
\end{definition}

We show that this notion of explicitness can be achieved by our 
construction. It will be convenient to give a similar definition for 
a signing of a graph, and for a composition of such functions:
\begin{definition}
Let $h_n:\{0,1\}^t \times [n] \times [n] \rightarrow \{-1,1\}$, with
$t = O(\log n)$. Given $r \in \{0,1\}^t$, and
a graph $G$ on $n$ vertices,
$h_n$ defines a signing $s_{h_n}$ of $G$ by 
$s_{h_n}(r)(i,j) = h_n(r,i,j)$.
We say that $h_n$ is
a {\em $\delta$-probabilistic strongly explicit} description of an 
$(n,d)$-good signing,
if given $n$, $h_n$ can be computed in time $\polylog(n)$, and, 
for any $(\log n, \gamma(d))$-sparse $d$-regular graph $G$ on $n$ vertices,
with probability at least $1-\delta$  (over
a uniform choice of $r$), $h_n$ defines an $(n,d)$-good signing.
\end{definition}
\begin{definition}
Let  $f_n:\{0,1\}^{t_1} \times [n] \times [n] \rightarrow \{0,1\}$,
and $h_n:\{0,1\}^{t_2} \times [n] \times [n] \rightarrow \{-1,1\}$ be as
above.
Their {\em composition}, 
$f_{2n}:\{0,1\}^t \times [2n] \times [2n] 
\rightarrow \{0,1\}$, with $t = \mbox{max}\{t_1,t_2\}$ is as follows. For $r \in \{0,1\}^t$, 
let $r_1$ be the first $t_1$ bits of $r$, and $r_2$ the first $t_2$ bits in $r$.
$f_{2n}$ is such that the graph $G_{f_{2n}}(r)$ is the $2$-lift
of $G_{f_n}(r_1)$ described by the signing $s_{h_n}(r_2)$. 
\end{definition}

The following lemma is easy, and we omit the proof:
\begin{lemma}\label{lem-compo}
Let $f_n$ be a $\delta_1$-probabilistic strongly explicit 
description of an $(n,d)$-good expander, and $h_n$ a 
$\delta_2$-probabilistic strongly 
explicit description of an $(n,d)$-good signing. Then their composition
is a  $(\delta_1+\delta_2)$-probabilistic strongly explicit description of 
an $(2n,d)$-good expander.
\end{lemma}

Think of an $(\epsilon,k)$-wise independent space $\Omega_m$ as a function
$\omega:\{0,1\}^t \rightarrow \{-1,1\}^m$, where $|\Omega_m| = 2^t$. It follows
from the work of Naor and Naor (\cite{Naor}), that not only can $\omega$ be
computed efficiently, but that given $r \in \{0,1\}^t, p \in [m]$
$\omega(r)_p$ (the $p$'th
coordinate of $\omega(r)$) can be computed efficiently (i.e. in time
$\polylog(m)$). Take $m = {n \choose 2}$, and think of the elements of 
$\{-1,1\}^m$ as being indexed by unordered pairs $(i,j) \in {[n] \choose 2}$.
Define $h_n(r,i,j) =  \omega(r)_{i,j}$. It follows from the above discussion
than $h_n$ is a $\frac 1 {n^2}$-probabilistic strongly explicit 
description of an $(n,d)$-good signing, for $k$ and $\epsilon$ as above.\\
We now describe how to construct a $\delta$-probabilistic strongly explicit 
description of an $(N,d)$-good expander graph.
Let $G$ be an $(n,d)$-good expander, with $n \geq \frac 1 {\delta}$. For 
$i=0,\ldots,l = log(N/n)$,
define $n_i = n \cdot 2^i$ and $m_i = {n_i \choose 2}$. 
Define $k_i = d \log n_i$. Let
$\omega_i:\{0,1\}^{t_i} \rightarrow \{-1,1\}^{m_i}$ be a description of
an $(\epsilon_i,k_i)$-wise independent space of bit strings of length $m_i$,
where $\epsilon_i$ is such that an element chosen uniformly at random
from this space yields an $(n_i,d)$-good signing with probability at least
$1-\frac 1 {n_i^2}$.\\
The functions  $h_{n_i}(r,p,q) =  \omega_i(r)_{p,q}$ are
$\frac 1 {n_i^2}$-probabilistic strongly explicit 
descriptions of $(n_i,d)$-good signings. Let $f_n$ be a description of $G$.
For simplicity, assume that adjacency in $G$ can be decided in time
$\polylog(n)$. Thus, $f_n$ is, trivially, a $0$-probabilistic strongly explicit 
description of an $(n,d)$-good expander. Define $f_{n_i}$ as the composition of
$f_{n_{i-1}}$ and $h_{n_{i-1}}$. It follows from this construction
and Lemma \ref{lem-compo} that:
\begin{lemma}\label{lem-pse}
$f_{n_l}$ is an $\frac 1 n$-probabilistic strongly explicit 
description of an $(N,d)$-good expander graph.
\end{lemma}

\section{A converse to the Expander Mixing Lemma}\label{sec-lme}
So far, we have concentrated on the algebraic approach to expansion
in graphs. Namely, that a graph is an
$(n,d,\lambda)$-expander if all eigenvalues but the largest are, in
absolute value, at most $\lambda$.
A combinatorial definition says that a $d$-regular
graph on $n$ vertices is an $(n,d,c)$-edge expander if every set 
of vertices, $W$, of size at most $n/2$,
has at least $c|W|$ neighbors emanating from it.\\
The two notions are closely related. Thus (cf. \cite{AlSp}),
an $(n,d,\lambda)$-expander is also an 
$(n,d,\frac{d-\lambda} {2})$-edge expander. 
Conversely, an  $(n,d,c)$-edge expander is also an
$(n,d,d-\frac {c^2} {2d})$-expander\footnote{A related result, showing 
that {\em vertex} expansion implies spectral gap appears in \cite{Al86}. The
implication from edge expansion is easier, 
and the proof we are aware of is also due to Noga Alon.}. 
Thus, though the two notions
of expansion are qualitatively equivalent, they are far from
being quantitatively the same. While algebraic expansion yields
good bounds on edge expansion, the reverse implications
are very weak. It is also known that this is not just a failure
of the proofs and indeed this estimate is nearly
tight \cite{Kahale}. Is there, we ask, another combinatorial property that
is equivalent to spectral gaps? We next answer this question.\\

For two subsets of vertices, $S$ and $T$, let $e(S,T)$ denote the number
of edges between them. We follow the terminology of \cite{Thomason}:
\begin{definition}\label{def-jumbled}
A $d$-regular graph $G$ on $n$ vertices is {\em $(d,\alpha)$-jumbled}, if
for every two subsets of vertices, $A$ and $B$,
\[|e(A,B) - d|A||B|/n| \leq \alpha \sqrt{|A||B|}.\]
\end{definition}

A very useful property of $(n,d,\lambda)$-expanders, known as the
Expander Mixing Lemma (cf. \cite{AlSp}), is that a  
an $(n,d,\lambda)$-expander is $(d,\lambda)$-jumbled.
Lemma \ref{wonder} implies the promised converse to this well known fact:

\begin{corollary}\label{lme}
Let $G$ be a $d$-regular graph on $n$ vertices. Suppose that for any
$S,T \subset V(G)$, with $S \cap T = \emptyset$
\begin{eqnarray*}
|e(S,T) - \frac {|S||T|d} {n}| \leq \alpha \sqrt{|S||T|}
\end{eqnarray*}
Then all but the largest eigenvalue of $G$ are bounded, in absolute
value, by $O(\alpha (1+\log(d/\alpha)))$. 
\end{corollary}
\begin{note}\label{note-jumbled}
In particular, this means that for a $d$-regular graph $G$,
$\lambda(G)$ is a $\log d$ approximation of the ``jumbleness'' parameter
of the graph.
\end{note}

\proof (Corollary \ref{lme})
Let $A$ be the adjacency matrix of $G$. Denote $B = A - \frac d n J$, where
$J$ is the all ones $n \times n$ matrix. Clearly $B$ is symmetric, and
the sum of the absolute value of the entries in each row is at most $2d$. 
The first eigenvalue of $A$ is $d$. The other eigenvalues  of $A$ are also 
eigenvalues  of $B$.
Thus, for the corollary to follow from Lemma \ref{wonder} it suffices to
show that for any two vectors, $u,v \in \{0,1\}^n$:
\[|uBv| = |u\frac d n J v - uAv| \leq \alpha ||u|| ||v||.\]
This is exactly the hypothesis for the sets $S(u)$ and $S(v)$.
\QED

The corollary is actually tight, up to a constant multiplicative factor, as
we now show:
\begin{theorem}\label{tight-lme}
For any large enough $d$, and $7 \sqrt{d} < \alpha < d$, there exist 
infinitely many $(d, \alpha)$-jumbled graphs with second eigenvalue 
$\Omega(\alpha (\log(d/\alpha) + 1))$.
\end{theorem}

It will be useful to extend Definition \ref{def-jumbled} to unbalanced 
bipartite graphs: 
\begin{definition}
A bipartite graph $G = (U,V,E)$ 
is {\em $(c,d,\alpha)$-jumbled}, if
the vertices in $U$ have degree $c$, those in $V$ have degree $d$, and
for every two subsets of vertices, $A \subset U$ and $B \subset V$,
\[|e(A,B) - d|A||B|/|U|| \leq \alpha \sqrt{|A||B|}.\]
\end{definition}

We note that such bipartite graphs exit:
\begin{lemma}
For $c | d$ and $\alpha = 2\sqrt{d}$, 
there exist $(c,d,\alpha)$-jumbled graphs.
\end{lemma}
\proof
Let $G'=(U',V',E')$ be a $c$-regular
Ramanujan bipartite graph, such that $|U'| = |V'| = n$. Let
$G = (U,V,E)$ be a bipartite graph obtained from $G'$ by partitioning
$V'$ into subsets of size $d/c$, and merging each subset into a vertex,
keeping all edges (so this is a multi-graph) .\\
Let $A \subset U$ and $B \subset V$. Let $A' = A$, and $B'$ the set of
whose merger gives $B$. Clearly $e(A,B) = e(A',B')$, $|B'| = d/c |B|$.
As $G'$ is Ramanujan,
by the expander mixing lemma
\[|e(A',B') - c|A'||B'|/n| \leq 2 \sqrt{c |A'||B'|},\]
or:
\[|e(A,B) - d|A||B|/|U|| \leq 2 \sqrt{d |A||B|}.\]
\QED

We shall need the following inequality, that can be easily proven by induction:
\begin{lemma}\label{lem-sum}
For $i=0,\ldots,t$ let $a_i$ be numbers in $[0,2^{2i}\cdot N]$, for some $N>0$.
Then
\[(\sum a_i 2^{-i})^2 \leq 3 N \sum a_i.\]
\end{lemma}

\proof (Theorem \ref{tight-lme})
Fix $d$, and $\sqrt{d} < \Delta < d$. 
Set $t = \half \log (\frac {3d} {4 \Delta})$, $\tau = \sum_{i=0}^t 4^i = 
\frac d \Delta < \Delta$. Let $N$ be some large number and $n = \tau N$ 
(this will be the number of vertices). For
$i,j=0,\ldots,t$ set $d_{i,j} = \frac d \tau 2^{2j} + \Delta
2^{j-i}$ when $i,j < t$ or $i=j=t$, and $d_{i,j} = \frac d
\tau 2^{2j} - \Delta 2^{j-i}$ otherwise. Note that in this case
$\frac d \tau 2^{2j} - \Delta 2^{j-i} \geq 2^j(\Delta 2^j - \Delta 2^{-i}) >0$.
Set $\alpha_{i,j} = 2 \sqrt{\min\{d_{i,j},d_{j,i}\}}$.\\
Let $|V_i|$ be subsets of size $4^i N$, 
and $G$ a graph on vertices $V = \cup_{i=0}^t V_i$ (hence, $|V|$ = n).
For $0 < i,j \leq t$ construct
a $(d_{j,i},d_{i,j},\alpha_{i,j})$-jumbled
graph between $V_i$ and $V_j$ (or $(d_{i,i},\alpha_{i,i})$-jumbled if
$i=j$).\\
The theorem follows from the following two lemmata.
\begin{lemma}
$G$ is $(d,7 \Delta)$-jumbled.
\end{lemma}
\proof
It is not hard to verify that $G$ is indeed $d$ regular.
Take $A,B \subset V$, and denote
their size by $a$ and $b$. Denote $A_i = A \cap V_i$,
$B_i = B \cap V_i$, and their size by $a_i$ and $b_i$.
We want to show that:
\[|e(A,B) - dab/n| \leq 3 \Delta \sqrt{ab}.\]
For simplicity we show that $e(A,B) \leq dab/n + 7 \Delta \sqrt{ab}$.
A similar argument bounds the number of edges from below.
From the construction,
$|e(A_i,B_j) - d_{i,j}|A_i||B_j|/|V_j|| \leq \alpha_{i,j}\sqrt{|A_i||B_j|}$, 
or:
\[e(A_i,B_j) \leq d_{i,j}a_i b_j/(4^jN) + \alpha_{i,j}\sqrt{a_i b_j}.\]
Summing up over $i,j = 0,\ldots,t$ we get:
\begin{eqnarray*}
e(A,B) &\leq& \sum d_{i,j}a_i b_j/(4^jN) + 
\alpha_{i,j} \sum \alpha_{i,j}\sqrt{a_i b_j} \\
&\leq& d / n \sum a_i b_j + \Delta/N \sum a_i b_j 2^{-(i+j)}
+ \sum \alpha_{i,j}\sqrt{a_i b_j}.
\end{eqnarray*}
$\sum a_i b_j = ab$, so it remains to bound the error term.
As $a_i, b_i \in [0,N \cdot 2^{2i}]$, by Lemma \ref{lem-sum}:
\begin{eqnarray*}
\Delta/N \sum a_i b_j 2^{-(i+j)} &=&
\Delta/N (\sum a_i 2^{-i})(\sum b_i 2^{-i})\\
 &\leq&
\Delta/N (\sqrt{3N\sum a_i})(\sqrt{3N\sum b_i})=
3 \Delta \sqrt{ab}
\end{eqnarray*}
so it remains to show that
$\sum \alpha_{i,j} \sqrt{a_i b_j} \leq 4 \Delta \sqrt{ab}$,
and hence it's enough to show that:
\begin{eqnarray*}
\sum 2 \sqrt{\frac d \tau a_i b_j 2^{2\max\{i,j\}}} + 
\sum 2 \sqrt{\Delta 2^{|j-i|} a_i b_j}
\leq 4 \Delta \sqrt{ab}.
\end{eqnarray*}
Indeed, as $i,j \leq t$, and $2^{2t} < \tau$,
\begin{eqnarray*}
\sum \sqrt{\frac d \tau a_i b_j 2^{2\max\{i,j\}}} < \sqrt{d} \sum \sqrt{a_i b_j}
< \Delta \sqrt{ab}.
\end{eqnarray*}
Similarly, $2^{|j-i|} \leq 2^t < \sqrt{\tau} < \sqrt{\Delta}$, and so:
\begin{eqnarray*}
\sum \sqrt{\Delta 2^{j-i} a_i b_j} < \Delta \sum \sqrt{a_i b_j}
= \Delta \sqrt{ab}.
\end{eqnarray*}
\QED

\begin{lemma}
$\lambda(G) \geq \Delta (t+1)$
\end{lemma}
\proof
Take $x \in \R^n$ to be $-2^{-t}$ on vertices in $V_t$, and $2^{-i}$ on
vertices in $V_i$,
for $i<t$. It is easy to verify that $x \perp \vec{1}$, and that
$||x||^2 = N \cdot (t+1)$. Let $M$ be the adjacency matrix if $G$. 
Since $\vec{1}$ is an eigenvector of $M$ corresponding to the largest eigenvalue,
by the variational characterization of eigenvalues, 
$\lambda(G) \geq \frac {xMx} {||x||^2}$.
Hence, to prove the lemma it suffices to show that $xMx \geq \Delta N (t+1)^2$. Indeed:
\begin{eqnarray*}
xMx &=& \sum_{i,j = 0}^t d_{i,j} 4^i N 2^{-(i+j)}
- 4 \sum_{i=0}^{t-1} d_{i,t} 4^i N 2^{-(i+t)}\\
&=& \frac d \tau N \sum_{i,j = 0}^t 2^{i+j}
+ \Delta N \sum_{i,j = 0}^t 1
- 4 \frac d \tau N \sum_{i=0}^{t-1} 2^{i+t}
+ 4 \Delta N \sum_{i=0}^{t-1} 1\\
&=& \frac d \tau N (2^{2(t+1)} - 4 \cdot 2^{2t}) + \Delta N ((t+1)^2 + 4t) 
> \Delta N (t+1)^2.
\end{eqnarray*}
\QED

\ignore{
It is known that for a random $d$-regular graph, w.h.p., the condition
in Corollary \ref{lme} holds with $\alpha = O(\sqrt{d})$ (cf. \cite{Shlomo}).
Hence, it follows from the corollary that w.h.p., such a graph is an 
$(n,d,O(\sqrt{d}\log{d}))$-expander. This result is weaker than 
previous ones (\cite{FKS,Fri91,Fri03}) by a multiplicative factor
of $\log d$.
}

\section{Reflections on Lemma \ref{wonder}}
\subsection{Finding the proof: LP-duality}
As the reader might have guessed, the 
proof for Lemma \ref{wonder} was discovered by
formulating the problem as a linear program.
Define
$\Delta_{i,j} = |x^i A x^j|$. Our assumptions translate to: 
\[\forall 1 \leq i \leq j \leq k:
|\Delta_{i,j}| \leq \alpha \sqrt{s_i s_j},\]
\[\forall 1 \leq i \leq k:
\sum_j|\Delta_{i,j}| \leq d s_i.\]

We want to deduce an upper bound on $|xAx|$. 
In other words, we are asking, 
under these constraints, how big 
\begin{eqnarray*}
\frac {|xAx|}  {||x||^2} \leq \frac {\sum_{i,j=1}^k \Delta_{i,j} 2^{-(i+j)}} 
{\sum_i 2^{-2i} s_i}
\end{eqnarray*}
can be. 

The dual program is to minimize:
\[\frac {\alpha \sum_{i \leq j} b_{i,j} \sqrt{s_i s_j} + d \sum_i c_i s_i} 
{\sum_i 2^{-2i} s_i}\]
under the constraints:
\begin{eqnarray*}
&\forall 1 \leq i < j \leq k,\;\;& b_{i,j} + c_i + c_j \geq 2^{-(i+j)+1}\\
&\forall 1 \leq i \leq k,& b_{i,i} + c_i \geq 2^{-2i}\\
&\forall 1 \leq i \leq j \leq k,& b_{i,j} \geq 0\\
&\forall 1 \leq i \leq k,& c_i  \geq 0\\
\end{eqnarray*}
The following choice of $b$'s and $c$'s satisfies the constraints, and gives 
the desired bound. These indeed appear in the proof above: 
\begin{eqnarray*}
&\forall 1 \leq i < j \leq k, j < i + \gamma,\;\;& b_{i,j} = 2^{-(i+j)+1}\\
&\forall 1 \leq i < j \leq k, j \geq i + \gamma,& b_{i,j} = 0\\
&\forall 1 \leq i \leq k,& b_i = 2^{-2i}\\
&\forall 1 \leq i \leq k,& c_i = 2^{-2i -\gamma + 1}\\
\end{eqnarray*}
Note that we can further constrain the $\delta$'s, 
and perhaps get better
bounds this way. 

\subsection{Algorithmic aspect}
Lemma \ref{wonder}
is algorithmic, in the sense that given a matrix with
a large eigenvalue, we can efficiently construct, from its eigenvector, a pair
$u, v \in \{0,1\}^n$ such that $S(u) \cap S(v) = \emptyset$, and 
$|uAv| \geq \alpha ||u|| ||v||$ (There is a small caveat - in the proof we
used a probabilistic argument for rounding the coordinates. 
This can be easily derandomized using the conditional probabilities method). 
Taking into consideration Note \ref{note-jumbled}, given a $d$-regular
graph $G$ where $\lambda(G)$ is large, one can efficiently find disjoint
subsets $S$ and $T$, such that $e(S,T) - \frac d n |S| |T| \geq 
c \cdot \lambda(G)/\log d \sqrt{|S||T|}$ (for some constant $c$).
We speculate that this might be useful in designing graph 
partitioning algorithms.

In Lemma \ref{lem-2} we showed the with positive probability a random signing
has a small spectral radius. It is interesting to find a sample space where 
this happens with high probability.
We conjecture that the following algorithm works: Choose a signing at random.
If some eigenvalue of the signed matrix is big,
use Lemma \ref{wonder} to find a pair $u,v \in \{0,1\}^n$ such that
$|uAv| \geq \alpha ||u|| ||v||$. 
Choose a signing at random for the edges between $S(u)$ and $S(v)$.
Repeat until all eigenvalues are small.

\section{Acknowledgments}
We thank L\'aszl\'o Lov\'asz for insightful discussions, and
Efrat Daom for help with computer simulations. We thank Eran Ofek
for suggesting that Corollary \ref{lme} might be used to bound the
second eigenvalue of random $d$-regular graphs, and Avi Wigderson
for Lemma \ref{lem-pse}.
We are grateful for helpful comments given to
us by Alex Samorodnitsky, Eyal Rozenman and Shlomo Hoory.

\pagebreak

\bibliography{bib}

\pagebreak
\setcounter{section}{0}
\def\thesection{\Alph{section}}
\def\thesubsection{\Alph{section}.\arabic{subsection}}

\end{document}